\documentclass{amsart}

\usepackage{graphicx}
\usepackage{amssymb}
\usepackage{enumerate}
\usepackage{amsmath}

\newtheorem{theorem}{Theorem}[section]
\newtheorem{lemma}[theorem]{Lemma}
\newtheorem{proposition}[theorem]{Proposition}

\theoremstyle{definition}
\newtheorem{definition}[theorem]{Definition}

\theoremstyle{remark}
\newtheorem{remark}[theorem]{Remark}

\numberwithin{equation}{section}

\begin{document}

\title[Geometric Lorenz attractors]{Parameter-shifted shadowing property for geometric Lorenz attractors}

\author{Shin Kiriki} 
\address{Department of Mathematical Sciences, Tokyo Denki University, Hatoyama, Hiki, Saitama-ken, 350-0394, JAPAN}
\email{ged@r.dendai.ac.jp}
\thanks{The first author was supported in part by 
Research Institute for Science and Technology at TDU Grant Q02J-02, Q03J-08.}

\author{Teruhiko Soma}
\address{Department of Mathematical Sciences, Tokyo Denki University, Hatoyama, Hiki, Saitama-ken, 350-0394, JAPAN}
\email{soma@r.dendai.ac.jp}

\date{April 10, 2003 and, in revised form, July 31, 2003}

\begin{abstract}
In this paper, we will show that any geometric Lorenz flow in a definite class satisfies the parameter-shifted shadowing property.
\end{abstract}

\maketitle

\section{Introduction}\label{S1}
We will study the problem whether there exists a definite class of geometric Lorenz flows which can be depicted as accurately as one desires.
Theoretically, such an accurate depiction is guaranteed by the shadowing property.
However, Komuro \cite{Ko} showed that geometric Lorenz flows do not satisfy the (parameter-fixed) shadowing property except very restricted cases.
So, we need to consider our problem under a somewhat relaxed condition, which is the parameter-shifted shadowing property in our case.

The geometric Lorenz model is one of important examples in  dynamical systems, 
which was studied in the initial stages by 
Guckenheimer and Williams \cite{G, W1, GW, W2}, Afraimovich-Bykov-Shil'nikov \cite{ABS} 
and Yoke-Yoke \cite{Y}.
Their aim was   to construct     topologically a simple  mechanism   
which can give results similar to that of the parametrized ODE system in $\mathbb{R}^3$ presented experimentally by Lorenz \cite{L}. 
For some parameter values, 
Lorenz observed  typical characters of chaotic motions in  butterfly-shaped   attractors. 
The question whether or not the original Lorenz system for such parameter values has the same structure as the  geometric Lorenz model has been unsolved for more than 30 years.  
By combination of normal form theory and rigorous computations, 
Tucker \cite{T} answered this question affirmatively,  that is, 
for classical parameters, 
the original Lorenz system has a robust strange attractor  which is given by the same rules as for the geometric Lorenz model.
From these facts, we know that the geometric Lorenz model is crucial in the study of 
Lorenz  dynamical systems.  
See Viana \cite{V} for more information.

The first return map on a Poincar\'e cross section of a geometric Lorenz flow is a Lorenz map $L:\Sigma\setminus \Gamma\longrightarrow \Sigma$, where $\Sigma=\big\{(x, y)\in \mathbb{R}^2;\ \vert x\vert, \vert y\vert \leq 1 \big\}$ and $\Gamma=\big\{(0, y)\in \mathbb{R}^2;\ \vert y\vert \leq 1 \big\}$.
So, we will first prove the parameter-shifted shadowing property (PSSP) for Lorenz maps.\\

\noindent{\bf Theorem A}
{\it 
There exists a definite set $\mathcal{L}$ of Lorenz maps satisfyingly the following condition.
\begin{itemize}
\item
For any $L\in \mathcal{L}$ and any $\varepsilon >0$, there exist $\mu >0$ and $\delta >0$ such that any $\delta$-pseudo orbit of the Lorenz map $L_\mu$ with $L_\mu (x,y)=L(x,y)-(\mu x,0)$ is $\varepsilon$-shadowed by an actual orbit of $L$.\end{itemize}
}

The strict description of $\mathcal{L}$ is given in the next section.

In the case when any elements in a one-parameter family $\{f_\mu\}_{\mu\in I}$ are naturally defined maps,
 `PSSP for $f=f_0$' means that any $\delta$-pseudo-orbit for $f$ is $\varepsilon$-shadowed 
 by an actual orbit of $f_\mu$ for some $\mu\in I$. This idea was first introduced by 
Coven-Kan-Yorke \cite{CKY} and 
Nusse-Yorke \cite{NY} in some one-dimensional dynamics. 
See also Kiriki-Soma \cite{KS} for PSSP for Lozi maps.
In the present case, $L_\mu$'s other than the original $L$ are artificially defined maps.
We wish here to describe actual orbits of the given map $L$ as accurately as possible but not those of auxiliary maps $L_\mu$, $\mu>0$.
Thus, we adopt as our definition of PSSP for $L$ that any $\delta$-pseudo-orbit for $L_\mu$ is $\varepsilon$-shadowed by an actual orbit of $L$.

As an application of Theorem A, we have the following result, which is our main theorem.\\

\noindent{\bf Theorem B}
{\it 
Any geometric Lorenz flow controlled by a Lorenz map $L\in \mathcal{L}$ has the parameter-shifted shadowing property.\\
}

See the next section for the definition of the parameter-shifted shadowing property of Lorenz flows.

\section{Preliminaries}

\begin{figure}[hbtp]
\centering
\includegraphics[bb=50 523 470 770,clip,width=9.3cm]{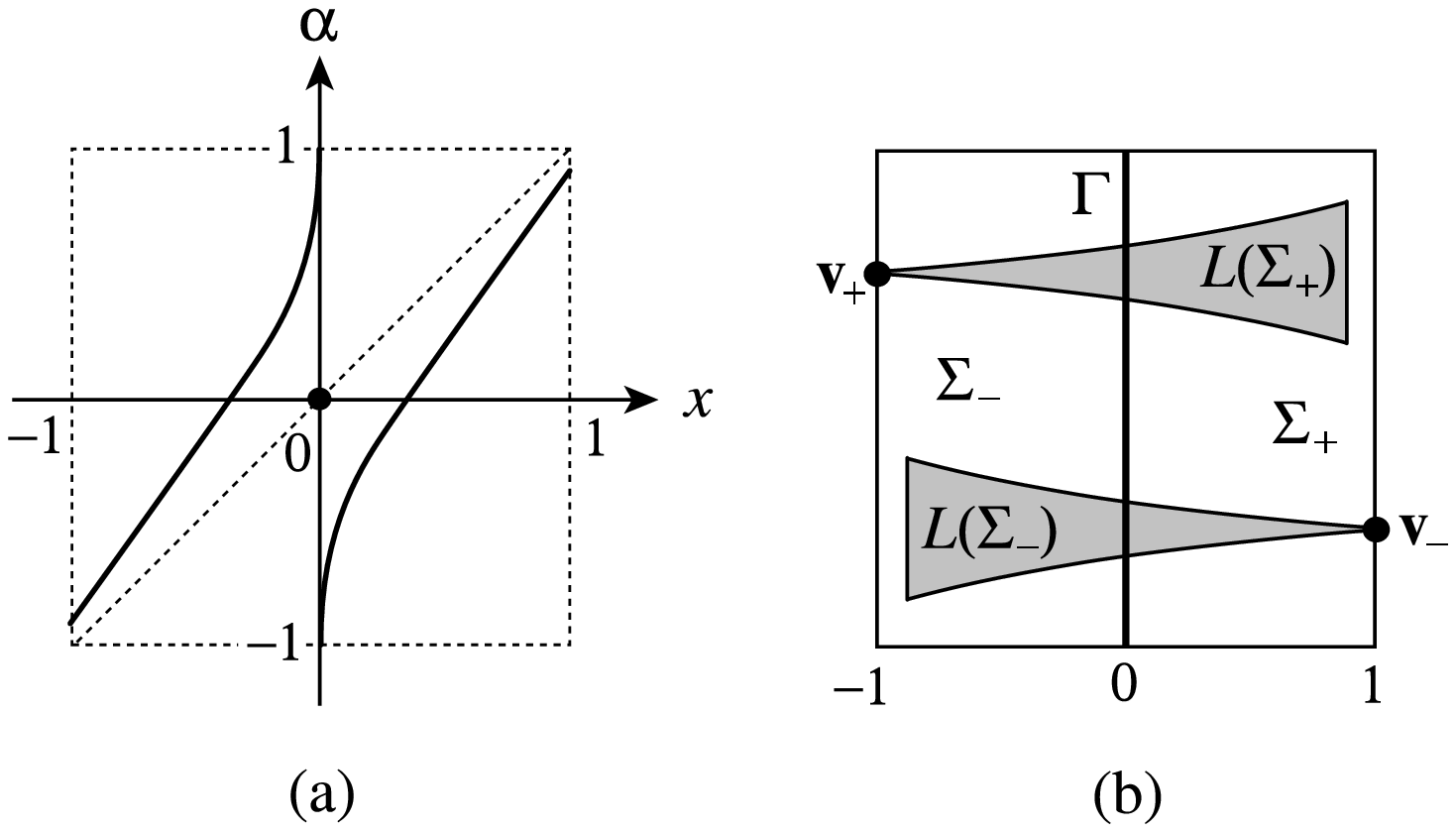}
\caption{}
\label{fig1}
\end{figure}

Let $\Sigma_\pm$ denote the components of $\Sigma\setminus \Gamma$ with $\Sigma_\pm\ni (\pm 1,0)$.
A map $L : \Sigma\setminus \Gamma\longrightarrow\Sigma$ 
is said to be a \textit{Lorenz map} if  
it is a  piecewise $C^1$   diffeomorphism which has  the following form 
$$L(x, y)=\big(\alpha(x), \beta(x,y)\big),$$ 
where  
$\alpha : [-1, 1]\setminus\{0\}\rightarrow  [-1, 1]$  is  a piecewise $C^{1}$-map  
with symmetric property $\alpha(-x)=-\alpha(x)$ and satisfying 
\begin{equation} \label{L-condition1}
\left\{\begin{array}{l}
	\lim_{x\rightarrow 0+}\alpha(x)=-1, \quad   \alpha(1)<1,\\
\lim_{x\rightarrow 0+}\alpha'(x)= \infty,\quad  \alpha'(x)>\sqrt{2}\ \mbox{for any}\ x\in (0, 1], 
\end{array}\right.
\end{equation}
see Fig.\ \ref{fig1}-(a), and 
$\beta:\Sigma\setminus \Gamma\longrightarrow [-1,1]$ is a contraction in the $y$-direction.
Moreover, it is required that the images $L(\Sigma_+)$, $L(\Sigma_-)$ are mutually disjoint cusps in $\Sigma$, where the vertices $\mathbf{v}_+,\mathbf{v}_-$ of $L(\Sigma_\pm)$ are contained in $\{\mp 1\}\times [-1,1]$ respectively, see Fig.\ \ref{fig1}-(b).

Now, we introduce the notion of the shadowing property for  Lorenz planar maps.

\begin{definition}{\rm
For $\delta>0$, 
a sequence $\{\mathbf{x}_n\}_{n\geq 0}\subset \Sigma$ is called 
a $\delta$-\textit{pseudo-orbit} of a Lorenz map $L$ if  
$$\vert L(\mathbf{x}_n)-\mathbf{x}_{n+1}\vert\leq \delta$$
 for any integer $n\geq 0$. 
Here, we suppose that  if $\mathbf{x}_n\in \Gamma$, then 
$\mathbf{x}_{n+1}$ is contained in one of the $\delta$-neighborhoods of $\mathbf{v}_+$ and $\mathbf{v}_-$.
}
\end{definition}

\begin{definition}{\rm 
A Lorenz map $L$  has the  \textit{parameter-shifted  shadowing property}, for short \textit{PSSP}, 
if 
there exists a one-parameter family $\{L_\mu\}_{\mu\in I}$ of Lorenz maps with $I=[0,\mu_0]$ for $0\leq \mu_0<1$ satisfying following conditions (i) and (ii).
\begin{enumerate}
	\renewcommand{\labelenumi}{(\roman{enumi})}
	\item $L_0=L$.
	\item For any $\varepsilon>0$, there exist  $\delta>0$ and $\mu\in I$ such that 
	any $\delta$-pseudo-orbit $\{\mathbf{x}_n\}_{n\geq 0}$ of $L_{\mu}$ 
	is \textit{$\varepsilon$-shadowed} by an actual orbit of $L$, i.e.\ 
	there exists a $\mathbf{z}\in \Sigma\setminus \Gamma $ such that 
	$$\vert L^n(\mathbf{z})-\mathbf{x}_{n}\vert\leq \varepsilon$$ for any $n\geq 0$.
\end{enumerate}
}
\end{definition}

When the parameter of $\{L_\mu \}_{\mu\in I}$ is fixed, i.e.\ $I=\{0\}$, 
the definition of PSSP is identical to that of the original (parameter-fixed) shadowing property given in \cite{An}.
According to Komuro \cite[Theorem 1]{Ko}, $L$ has the parameter-fixed shadowing property only when $\alpha(1)=1$, 
see also \cite{Pi}.
In our case, since $\alpha(1)<1$ by (\ref{L-condition1}), any Lorenz map $L$ does not have the original shadowing property.

We are mainly concerned with Lorenz maps $L(x,y)=(\alpha(x),\beta(x,y))$ satisfying the following extra conditions (\ref{L-condition2})-(\ref{L-condition4}).
\begin{equation} \label{L-condition2}
  \left|  \frac{\partial \beta}{\partial x}(x,y)\right|, \ \left|  \frac{\partial \beta}{\partial y}(x,y)\right|<\frac{3}{4\sqrt{2}}\quad \mbox{for any}\ (x,y)\in \Sigma\setminus\Gamma,
\end{equation}
\begin{equation}\label{L-condition3}
0.8<\alpha^2(1)<\alpha(1)<1,
\end{equation}
\begin{equation}\label{L-condition4}
\ \alpha'(x)<2\quad \mbox{for any}\ x \mbox{ with }  0.8< x\leq 1.
\end{equation}
These conditions are not so severe, and it is not hard for us to construct various Lorenz maps satisfying them practically.
In the condition (\ref{L-condition2}), we took the concrete value $3/(4\sqrt{2})$ in order to simplify the proof of the theorem below.
In fact, one can prove the theorem under the weaker assumption:
$$\sup_{(x,y)\in \Sigma\setminus \Gamma}\left\{ \left|  \frac{\partial \beta}{\partial x}(x,y)\right|,\ \left|  \frac{\partial \beta}{\partial y}(x,y)\right|\right\}<\frac{1}{\sqrt{2}}.$$

The following is the precise statement of Theorem A.

\begin{theorem}[PSSP for  Lorenz planar maps]\label{Th1}
Any Lorenz map $L$ with the conditions {\rm (\ref{L-condition2})-(\ref{L-condition4})} admits a one-parameter family $\{L_\mu\}_{\mu\in I}$,
$$L_\mu (x,y)=L(x,y)-(\mu x,0),$$
satisfying the parameter-shifted shadowing property.
Precisely, for any $\varepsilon>0$, there exist $\delta>0$ and $\mu\in I$ so that the following {\rm (i)} and {\rm (ii)} hold.
\begin{enumerate}
\renewcommand{\labelenumi}{\rm (\roman{enumi})}
\item
Any infinite $\delta$-pseudo-orbit $\{\mathbf{x}_n\}_{n=1}^\infty$ of $L_\mu$ {\rm (}possibly $\mathbf{x}_n\in \Gamma${\rm )} is $\varepsilon$-shadowed by the actual orbit $\{L^n(\mathbf{z})\}_{n=1}^\infty$ of $L$ for some $\mathbf{z}\in \Sigma$ with $\bigcup_{n=0}^\infty \{L^n(\mathbf{z})\}\cap \Gamma =\emptyset$.
\item
Any finite $\delta$-pseudo-orbit $\{\mathbf{x}_n\}_{n=1}^m$ of $L_\mu$ with $\mathbf{x}_m\in \Gamma$ is $\varepsilon$-shadowed by the actual orbit $\{L^n(\mathbf{z})\}_{n=1}^m$ of $L$ for some $\mathbf{z}\in \Sigma$ with $\bigcup_{n=0}^{m-1} \{L^n(\mathbf{z})\}\cap \Gamma =\emptyset$ and $L^m(\mathbf{z})\in \Gamma$.
\end{enumerate}
\end{theorem}

Let us identify $\Sigma$ with $\big\{(x, y, 1)\in \mathbb{R}^3;\ \vert x\vert, \vert y\vert \leq 1 \big\}$, 
 and $\Gamma$ with $\big\{(0, y, 1)\in \mathbb{R}^3;\ \vert y\vert \leq 1 \big\}$.
A $C^1$-vector field  $X_L$ on $\mathbb{R}^3$ is said to be a {\it geometric Lorenz vector field controlled by} a Lorenz map $L:\Sigma\setminus \Gamma\longrightarrow \Sigma$ if it satisfies the following conditions (i) and (ii).
\begin{enumerate}
\renewcommand{\labelenumi}{\rm (\roman{enumi})}
\item
For any point $(x,y,z)$ in a neighborhood of the origin $\mathbf{0}$ of $\mathbb{R}^3$, $X_L$ is given  by 
$(\dot{x}, \dot{y}, \dot{z})=(\lambda_1 x, -\lambda_2 y, -\lambda_3 z)$, 
where $\lambda_i$ are positive numbers  satisfying $\lambda_3<\lambda_1<\lambda_2$. 
Moreover, $\Gamma$ is contained in the stable manifold $W^s(\mathbf{0})$ of $\mathbf{0}$.
\item
All forward orbits of $X$ starting from $\Sigma\setminus\Gamma$ will return  to $\Sigma$ and the first return map is $L$.
\end{enumerate}

\begin{figure}[hbt]
\begin{center}
\includegraphics{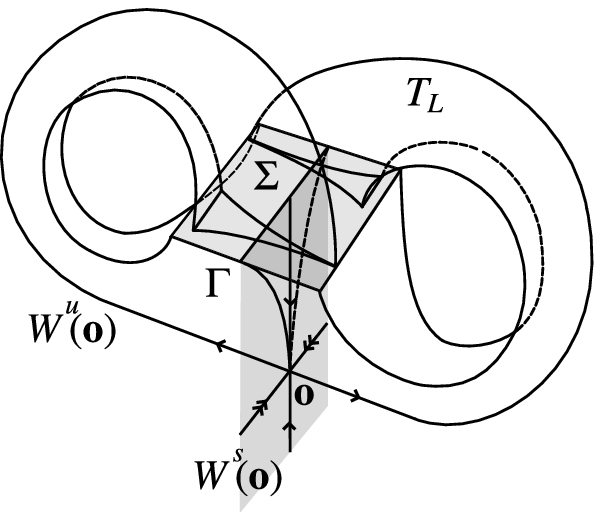}
\caption{}
\label{fig2}
\end{center}
\end{figure}

Note then that $\mathbf{0}$ is a singular point (an equilibrium) of saddle type, the local unstable manifold of $\mathbf{0}$ is tangent to  the  $x$-axis, and 
the local stable manifold of $\mathbf{0}$ is tangent to the  $yz$-plane as shown in Fig.\ \ref{fig2}. 
A $C^1$-map $\varphi_L:  \mathbb{R}^3\times \mathbb{R}\rightarrow \mathbb{R}^3$ is  the \textit{geometric Lorenz flow controlled by $L$} (for short $L$-\textit{Lorenz flow}) if it generated by 
  $X_L$, i.e.\ $\varphi_L(\mathbf{x}, 0)=\mathbf{x}$ and 
$(\partial/\partial t)\varphi_L(\mathbf{x}, t)=X_L(\varphi_L(\mathbf{x}, t))$.
The closure of $\bigcup_{\mathbf{z}\in \Sigma\setminus \Gamma}\varphi_L(\mathbf{z},[0,\infty))$ in $\mathbb{R}^3$ is homeomorphic to the genus two solid handlebody as illustrated in Fig.\ \ref{fig2}, which is called a \textit{trapping region} of $\varphi_L$ and denoted by $T_{\varphi_L}$ or $T_L$.
Any forward orbit for $\varphi_L$ with an initial point in $T_L$  can not escape from $T_L$.
The  invariant set $\bigcap_{t\geq 0} \varphi_L(T_L,t)$ for $X_L$ does not have any continuous hyperbolic splitting at  $\mathbf{0}$, but it 
belongs to an essential class called singular hyperbolic, 
which is  
studied  extensively  from various approaches  
 by Morales,  Pacifico and others, see  for details \cite{CM, CMP,  M, MP, MPP, MPP2} .

Now, we introduce the notion of PSSP for Lorenz flows.

\begin{definition}\label{flowcase}
Let $\psi$ be a geometric Lorenz flow, and $\delta,\tau$ positive numbers.
{\rm
\begin{enumerate}
\renewcommand{\labelenumi}{(\roman{enumi})}
\item
A sequence $\{\mathbf{x}_n\}_{n\geq 0}$ in $T_\psi$ with $\mathbf{x}_0\in \Sigma$ is a 
\textit{$(\delta,\tau)$-pseudo-orbit} for the flow  $\psi$ if 
there exists a sequence  $\{\tau_n\}_{n\geq 0}$ such that, for any $n\geq 0$,
$$
\tau\leq \tau_n\leq 2\tau\quad\mbox{and}\quad \left\vert \psi(\mathbf{x}_n,\tau_n)-\mathbf{x}_{n+1}\right\vert\leq \delta.
$$ 
For each $n\geq 0$, we set
$$\Psi_n=  \psi(\mathbf{x}_n, [0, \tau_n]) $$
and call $\{\Psi_n\}_{n\geq 0}$ the \textit{$(\delta,\tau)$-chain} for $\psi$ associated to $\{\mathbf{x}_n\}_{n\geq 0}$ (or more strictly to $\{\mathbf{x}_n;\tau_n\}_{n\geq 0}\}$).
\item
The $(\delta,\tau)$-chain $\{\Psi_n\}_{n\geq 0}$ is said to be \textit{$\varepsilon$-shadowed} by a flow 
$\varphi$
if there exists a point 
$ \mathbf{y}\in \Sigma$ and a surjective $C^1$-diffeomorphism $h:[0,\ \infty)\longrightarrow [0,\infty)$ satisfying
$$
	\left\vert
	\varphi(\mathbf{y}, h(t))-\psi(\mathbf{x}_n,t-\sum_{i=0}^{n-1}\tau_i)\right\vert\leq \varepsilon
$$
for any $t\geq 0$ with $\sum_{i=0}^{n-1}\tau_i\leq t\leq \sum_{i=0}^n\tau_i$.
Then, we also say that $\{\Psi_n\}_{n\geq 0}$ is $\varepsilon$-shadowed by $\varphi$ \textit{with} $\varphi(\mathbf{y},t); t\geq 0$.
\end{enumerate}
}
\end{definition}

\begin{figure}[hbtp]
\centering
\includegraphics[bb=73 583 512 798,clip,width=9.8cm]{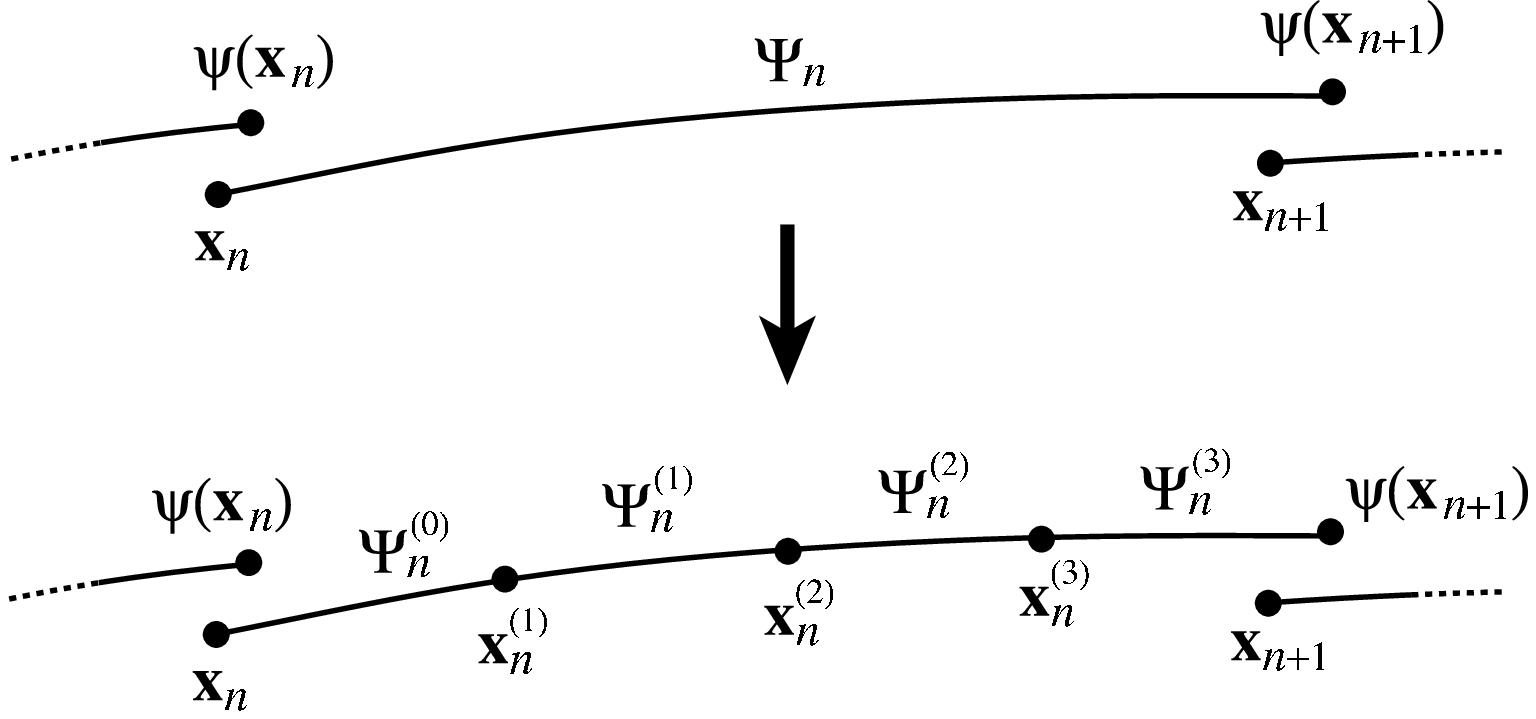}
\caption{}
\label{fig25}
\end{figure}

\begin{remark}{\rm
In Definition \ref{flowcase}\,(i), the upper bound condition $\tau_n\leq 2\tau$ is not essential, but added for our convenience.
When $\tau_n> 2\tau$ for some $n$, split $[0,\tau_n]$ into subintervals $[\tau_n^{(i-1)},\tau_n^{(i)}]$ $(i=1,\cdots,k)$ with $\tau_n^{(0)}=0$, $\tau_n^{(k)}=\tau_n$ and $\tau\leq \tau_n^{(i)}-\tau_n^{(i-1)}\leq 2\tau$.
Then, the expanded sequence of $\{\mathbf{x}_n\}_{n\geq 0}$ obtained by adding the entries $\mathbf{x}_n^{(i)}=\psi(\mathbf{x}_n,\tau_n^{(i)})$ $(i=1,\cdots,k-1)$ between $\mathbf{x}_n$ and $\mathbf{x}_{n+1}$ defines a $(\delta,\tau)$-pseudo-orbit for $\psi$ in the sense of Definition \ref{flowcase}\,(i), see Fig.\ \ref{fig25}.
}
\end{remark}

\begin{definition}{\rm
We say that 
a given geometric 
Lorenz flow $\varphi$ (or the vector field generating $\varphi$) has the \textit{parameter-shifted shadowing property} 
if  
there exists a $C^1$-one-parameter family $\{\varphi_\mu\}_{\mu \in [0,\mu_0]}$  of geometric 
Lorenz flows such that 
\begin{enumerate}
	\renewcommand{\labelenumi}{(\roman{enumi})}
	\item $\varphi_0=\varphi$;
	\item for any $\varepsilon>0$,  
	there exist $\delta,\tau>0$ and $\mu\in [0,\mu_0]$ such that  
	any $(\delta,\tau)$-chain for  $\varphi_{\mu}$  is $\varepsilon$-shadowed by $\varphi_0$.
\end{enumerate}
}
\end{definition}

The following is the precise statement of Theorem B.

\begin{theorem}[PSSP for Lorenz flows]\label{Th2}
Any geometric Lorenz flow controlled by a Lorenz map satisfying the conditions {\rm (\ref{L-condition2})-(\ref{L-condition4})} has the parameter-shifted shadowing property.
\end{theorem}

\begin{remark}[Absence of strong PSSP for Lorenz flows]\label{nonstrong}
{\rm 
Our PSSP for Lorenz flows is the weak one in the sense of Definition 3 in \cite{Ko}.
We say that a $(\delta,\tau)$-chain $\{\Phi_{\mu;n}\}_{n\geq 0}$ for $\varphi_{L_\mu}$ is {\it strongly} $\varepsilon$-shadowed by $\varphi_L$ with $\varphi_L(\mathbf{y},h(t))_{t\geq 0}$ if a diffeomorphism $h:[0,\infty)\longrightarrow [0,\infty)$ as in Definition \ref{flowcase} satisfies the extra condition: $\vert h^\prime(t)-1\vert<\varepsilon$ for any $t\geq 0$.
However, \textit{any} $\varphi_L$ as in Theorem \ref{Th2} has a constant $\varepsilon=\varepsilon(L)>0$ such that, for any $\delta,\tau>0$ and any $\mu \in I$ (possibly $\mu=0$), there exists a $(\delta,\tau)$-chain for $\varphi_{L_\mu}$ which is not strongly $\varepsilon$-shadowed by any actual flow of $\varphi_L$.
This implies that $\varphi_L$ does not have the strong PSSP.
In fact, one can define a $(\delta,\tau)$-pseudo-orbit $\{\mathbf{x}_n\}_{n\geq 0}$ in $T_{L_\mu}$ for $\varphi_{L_\mu}$ such that, in a small neighborhood of $\mathbf{0}$ in $\mathbb{R}^3$, the sequence satisfies $\mathbf{x}_{n_0}=\mathbf{x}_{n_0+1}=\cdots=\mathbf{x}_{n_0+m}$ for an arbitrarily large $m\geq 0$.
Such a sequence $\{\mathbf{x}_n\}_{n\geq 0}$ is not strongly $\varepsilon$-shadowed by $\varphi_L$.
The proof is elementary but somewhat tedious, so we will omit it.
}
\end{remark}

\section{PSSP  for  Lorenz planar maps}
Let $L:\Sigma\setminus \Gamma \longrightarrow \Sigma$ with $L(x,y)=(\alpha (x),\beta (x,y))$ be a Lorenz map satisfying the conditions (\ref{L-condition1})-(\ref{L-condition4}).
For any $\mu >0$, consider the function $\alpha_\mu:[-1,1]\setminus \{0\}\longrightarrow \mathbb{R}$ defined by 
$$\alpha_\mu(x)=\alpha(x)-\mu x.$$
If $\mu_0=\mu_0(\alpha)>0$ is sufficiently small, then for any $\mu\in [0,\mu_0]$, $\alpha_\mu$ is a function with $\alpha_\mu ([-1,1]\setminus \{0\})\subset [-1,1]$ and satisfies (\ref{L-condition1}), (\ref{L-condition3}) and (\ref{L-condition4}).
Set $I=[0,\mu_0]$.
Then, we have the one-parameter family $\{L_\mu\}_{\mu\in I}$ of Lorenz maps with 
$$L_\mu(x, y)=\big(\alpha_\mu(x),  \beta(x,y)\big)$$
for $(x, y)\in \Sigma\setminus\Gamma$.

By the condition (\ref{L-condition3}), there exists $0<\eta_0<1$ such that, for any interval $J=[-\eta,0]$ or 
$[0, \eta]$ with $0<\eta\leq \eta_0$ and any $\mu\in [0,\mu_0]$,
\begin{equation}\label{J}
\bigcup_{i=1}^3\alpha^i_\mu(J)\subset [0.8,1].
\end{equation}

For any $\varepsilon>0$, we set
\begin{equation}\label{epsilondelta}
\varepsilon_1=\min\left\{3\mu_0,\frac{\eta_0}{8},\frac{\varepsilon}{64}\right\}\quad\mbox{and}\quad \delta=\frac{\varepsilon_1}{100}.
\end{equation}

\begin{proposition}\label{alpha}
The map $\widehat\alpha=\alpha_{\varepsilon_1/3}$ satisfies the following {\rm (i)} and {\rm (ii)}.
\begin{enumerate}
\renewcommand{\labelenumi}{\rm (\roman{enumi})}
\item
Any infinite $\delta$-pseudo-orbit $\{x_n\}_{n=0}^\infty$ of $\widehat\alpha$ is $\varepsilon/8$-shadowed by an actual orbit $\{\alpha^n(z)\}_{n=0}^\infty$ of $\alpha$ for some $z\in [-1,1]$ with $\bigcup_{n=0}^\infty \{\alpha^n(z)\}\not\ni 0$.
\item
Any finite $\delta$-pseudo-orbit $\{x_n\}_{n=0}^m$ of $\widehat\alpha$ with $x_m=0$ is $\varepsilon/8$-shadowed by an actual orbit $\{\alpha^n(z)\}_{n=0}^m$ of $\alpha$ for some $z\in [-1,1]$ with $\bigcup_{n=0}^{m-1} \{\alpha^n(z)\}\not\ni 0$ and $\alpha^m(z)=0$.
\end{enumerate}
\end{proposition}

Consider any infinite $\delta$-pseudo-orbit $\{x_n\}_{n=0}^\infty$ of $\widehat\alpha$.
Let $l_0$ be the closed interval in $\mathbb{R}$ with $o(l_0)=x_0$ and $|l_0|=2\varepsilon_1$, where $|l_0|$ is the length of $l_0$ and $o(l_0)$ is the center of $l_0$.  
For the proof of Proposition \ref{alpha}, we will define a certain sequence of closed intervals $\{l_n\}_{n\geq 0}$ in $[-1,1]$ with $\mathrm{Int} l_n\cap \{0\}=\emptyset$ and $\alpha(l_n)\supset l_{n+1}$, 
where the notation $\mathrm{Int} l_n$ means the interior of $l_n$.

\begin{figure}[hbtp]
\centering
\includegraphics[bb=153 675 474 803,clip,width=7.1cm]{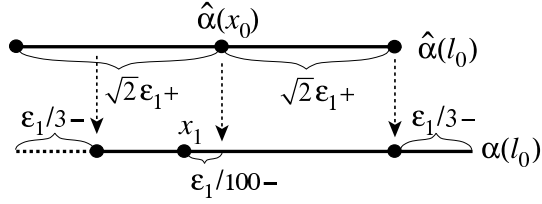}
\caption{`$\sqrt{2}\varepsilon_1+$' (resp.\ `$\varepsilon_1/3-$') in the figure means that the length between the corresponding points is at least $\sqrt{2}\varepsilon_1$ (resp.\ at most $\varepsilon_1/3$).
These rules are applied in any figures below.
}
\label{fig3}
\end{figure}

Suppose first that $l_0\cap \{0\}=\emptyset$.
We may assume that $l_0\subset (0,1]$.
Then, $\widehat\alpha(l_0)$ is the closed interval in $[-1,1]$ such that $\widehat\alpha(x_0)$ divides $\widehat\alpha(l_0)$ into two intervals of length at least $\sqrt{2}\varepsilon_1$.
Since $0\leq \alpha(x)-\widehat\alpha(x)=(\varepsilon_1/3)x\leq \varepsilon_1/3$ for any $x\in (0,1]$, $\alpha(l_0)$ is obtained from $\widehat \alpha(l_0)$ by $(0,\varepsilon_1/3)$-RHS-shifting, see Fig\ \ref{fig3},  
where  we say that
for two closed intervals $l=[e,f],l'=[e',f']$ and $0\leq \gamma\leq \eta$,  $l'$ is obtained from $l$ by $(\gamma,\eta)$-{\it right-hand-side shifting} (for short $(\gamma,\eta)$-RHS-shifting) if $\gamma \leq e'-e\leq \eta$ and $\gamma \leq f'-f\leq \eta$.

Since $|x_1-\widehat\alpha(x_0)|\leq \varepsilon_1/100$, 
the distance between $x_1$ and either end point of $\alpha(l_0)$ is at least $\sqrt{2}\varepsilon_1-\varepsilon_1/3-\varepsilon_1/100>\varepsilon_1$.
Thus, the interval $l_1$ with $|l_1|=2\varepsilon_1$ and $o(l_1)=x_1$ is contained in the {\it interior} of $\alpha(l_0)$.
If $l_1\cap \{0\}=\emptyset$, one can define $l_2\subset \mathrm{Int} \alpha(l_1)$ with $|l_2|=2\varepsilon_1$ and $o(l_2)=x_2$ similarly.

\begin{figure}[hbtp]
\centering
\includegraphics[bb=47 613 546 772,clip,width=11.0cm]{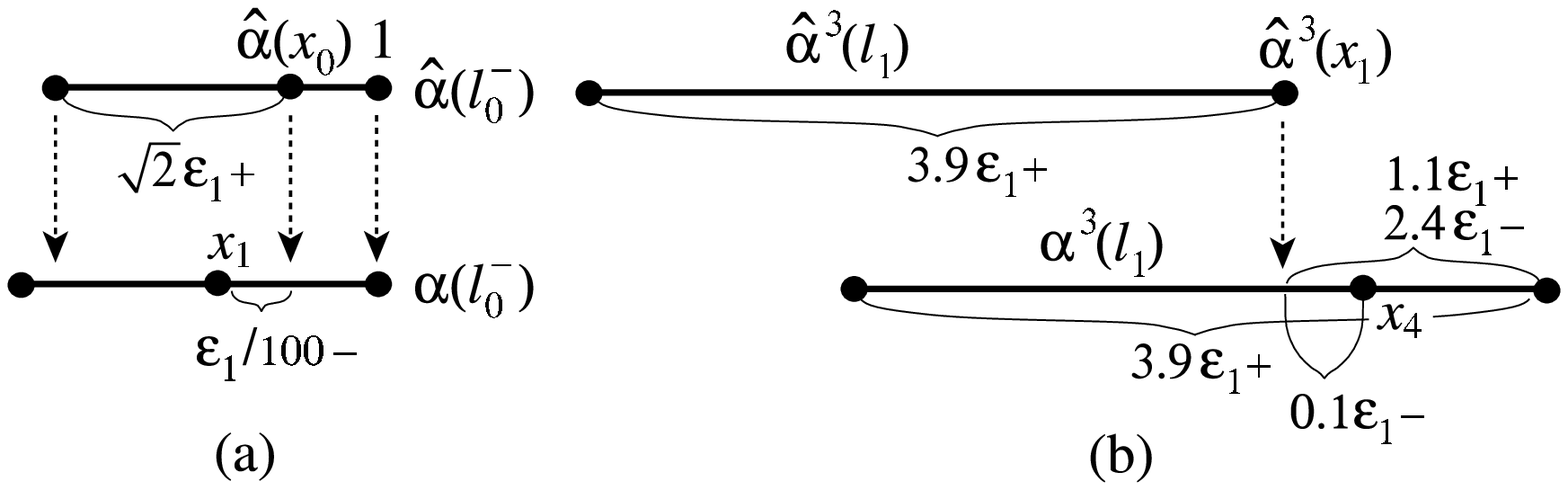}
\caption{}
\label{fig4}
\end{figure}

Suppose next that $l_0\cap \{0\}\neq \emptyset$.
We may assume that $x_0\leq 0$ and $\widehat\alpha(x_0)>0$.
Set $l_0^-=l_0\cap [-1,0]$.
Then, $\widehat\alpha(l_0^-)$ is the closed interval in $[0.8,1]$ containing $1$, see Fig.\ \ref{fig4}-(a).
Note that the distance between $\widehat\alpha(x_0)$ and the end point of $\widehat\alpha(l_0^{-})$ other than $1$ is at least $\sqrt{2}\varepsilon_1$.
The interval $\widehat\alpha(l_0^-)$ is obtained from $\alpha(l_0^-)$ by $(0,\varepsilon_1/3)$-RHS-shifting.
Since $|\widehat\alpha(x_0)-x_1|\leq \varepsilon_1/100$, the interval $l_1=[x_1-(\sqrt{2}-1/100)\varepsilon_1,x_1]$ is contained in $\alpha(l_0^-)$.
Note that, from the condition (\ref{J}), $\bigcup_{i=0}^2(\alpha^i(l_1)\cup\widehat\alpha^i(l_1))$ is contained in $[0.8,1]$.
Since $\alpha'(x)>\sqrt{2}$ for any $x\in (0,1]$, the length of the interval $\alpha^3(l_1)$ is at least $2\sqrt{2}(\sqrt{2}-1/100)\varepsilon_1>3.9\varepsilon_1$.
Since $0.8(2+\sqrt{2}+1)\varepsilon_1/3>1.1\varepsilon_1$ and $(2^2+2+1)\varepsilon_1/3<2.4\varepsilon_1$, $\alpha^3(l_1)$ is obtained from $\widehat\alpha^3(l_1)$ by $(1.1\varepsilon_1,2.4\varepsilon_1)$-RHS shifting, see Fig.\ \ref{fig4}-(b).
Since $\alpha'(x)<2$ for any $x\in [0.8,1]$,
$$|\widehat\alpha^3(x_1)-x_4|\leq (2^2+2+1)\frac{\varepsilon_1}{100}<0.1\varepsilon_1.$$
Thus, the closed interval $l_4$ with $|l_4|=2\varepsilon_1$ and $o(l_4)=x_4$ is contained in $\mathrm{Int} \alpha^3(l_1)$.
Set $l_2=\alpha(l_1)$ and $l_3=\alpha(l_2)$.
Then, $|l_2|<|l_3|\leq 2^2|l_1|<6\varepsilon_1$ and, for $i=1,2$, 
\begin{eqnarray*}
|\alpha^{i}(x_1)-x_{i+1}|&\leq& |\alpha^i(x_1)-\widehat\alpha^i(x_1)|+|\widehat\alpha^{i}(x_1)-x_{i+1}|\\
&\leq& (2+1)\left(\frac{\varepsilon_1}{3}+\frac{\varepsilon_1}{100}\right)<1.1\varepsilon_1.
\end{eqnarray*}
In particular, for any $y\in l_i$, $|y-x_i|\leq 6\varepsilon_1+1.1\varepsilon_1<8\varepsilon_1$. 
Now, for any given $\{x_n\}_{n=0}^\infty$ of $\widehat\alpha$, 
we get the sequence of closed intervals  $\{l_n\}_{n\geq 0}$ with $\mathrm{Int} l_n\cap \{0\}=\emptyset$ and $\alpha(l_n)\supset l_{n+1}$, 
For any $n\geq 0$, we set $l^{(n)}_{n+1}=\alpha^{-1}(l_{n+1})\cap l_n$.
For any $m>n+1$, $l_m^{(n)}$ can be defined inductively on $m-n$ by $l_m^{(n)}=\alpha^{-1}(l_m^{(n+1)})\cap l_n$.
Note that the restriction $\alpha^{m-n}|l_m^{(n)}:l_m^{(n)}\longrightarrow l_m$ is a bijection.

Lemma \ref{interval} is proved by applying the argument above repeatedly.

\begin{lemma}\label{interval}
There exists a sequence $\{l_n\}_{n=0}^\infty$ of closed intervals satisfying the following conditions {\rm (i)}-{\rm (iv)}.
\begin{enumerate}
\renewcommand{\labelenumi}{\rm (\roman{enumi})}
\item
$\varepsilon_1\leq |l_n|\leq 6\varepsilon_1$.
\item
For any $y\in l_n$, $|y-x_n|\leq 8\varepsilon_1$.
\item
If $l_n$ is not contained in $[-1,-0.8]\cup [0.8,1]$, then $|l_n|=2\varepsilon_1$ and $o(l_n)=x_n$.
\item
For any $n\geq 0$, $\alpha(l_n)$ contains $l_{n+1}$ and there exists $m>n$ with $l_m^{(n)}\subset \mathrm{Int} l_n$.
Moreover, if $l_n\ni 0$, then $\mathrm{Int} l^{(n)}_{n+1}\cap \{0\}=\emptyset$.
\end{enumerate}
\end{lemma}

The proof of Proposition \ref{alpha} follows easily from Lemma \ref{interval}.\\

\noindent\textit{Proof of Proposition \ref{alpha}.}
(i) Since $l_0\supset l_1^{(0)}\supset l_2^{(0)}\supset \cdots$, the intersection $\bigcap_{n=1}^\infty l_n^{(0)}$ is non-empty.
Take a point $z\in \bigcap_{n=1}^\infty l_n^{(0)}$.
Since $\alpha^n(z)\in l_n$ for any $n\geq 0$, by Lemma \ref{interval},
$$|\alpha^n(z)-x_n|\leq 8\varepsilon_1\leq \varepsilon/8.$$
If $0\in l_n$, then $\alpha^n(z)\in l_m^{(n)}\subset \mathrm{Int} l_{n+1}^{(n)}\subset l_n$ for some $m>n+1$.
This implies $\alpha^n(z)\neq 0$ for any $n\geq 0$.

(ii) Supposing that $\{x_n\}_{n=1}^m$ is a finite subsequence of an infinite $\delta$-pseudo-orbit  of $\widehat\alpha$, we have a sequence $l_0,l_1,\cdots,l_m$ of closed intervals satisfying the conditions (i)-(iv) of Lemma \ref{interval}.
If $x_m=0$, then $l_m$ is the interval with $|l_m|=2\varepsilon_1$ and $o(l_m)=0$.
In particular, $0\in \mathrm{Int} l_m$.
Thus, $\bigcap_{n=1}^m l_m^{(0)}$ contains a unique point $z$ with $\alpha^n(z)\neq 0$ for $0\leq n\leq m-1$ and $\alpha^m(z)=0$.
This completes the proof.
\hspace*{\fill}$\square$\\

\noindent\textit{Proof of Theorem \ref{Th1}.}
(i) Let $\{\mathbf{x}_n\}_{n=1}^\infty$ be any $\delta$-pseudo-orbit for $L_\mu$ with $\mu=\varepsilon_1/3$.
Set $[\mathbf{x}_n]_x=x_n$ and $[\mathbf{x}_n]_y=y_n$, where $[\mathbf{w}]_x$, $[\mathbf{w}]_y$ denote respectively the $x$ and $y$-coordinates of a point $\mathbf{w}\in \Sigma$.
Then, $\{x_n\}_{n=1}^\infty$ is a $\delta$-pseudo-orbit for $\widehat\alpha$.
By Proposition \ref{alpha}(i), there exists $z\in [-1,1]$ such that $\{\alpha^n(z)\}_{n=1}^\infty$ $\varepsilon/8$-shadows $\{x_n\}_{n=1}^\infty$.
If we set $\mathbf{z}=(z,y_0)$, then
$$|\mathbf{z}-\mathbf{x}_0|=|z-x_0|<\varepsilon.$$
Suppose that $|L^n(\mathbf{z})-\mathbf{x}_n|<\varepsilon$ for $n=0,1,\cdots,m$.
Let $J$ be a straight segment in $\Sigma$ connecting $L^m(\mathbf{z})$ with $\mathbf{x}_m$.
If $c:[0,\nu]\longrightarrow J$ be an arc-length parametrization of $J$, then $|\dot c|^2=\dot c_1^2+\dot c_2^2=1$ and $\nu<\varepsilon$, where $\dot c=(d/dt)c$ and $c(t)=(c_1(t),c_2(t))$.
Since $|\alpha^m(z)-x_{m}|<8\varepsilon_1\leq \eta_0$, $\mathrm{Int} J$ is disjoint from $\Gamma$.
In fact, if $\mathrm{Int} J\cap \Gamma$ were not empty, by (\ref{J}), 
then $|\alpha^{m+1}(z)-\alpha(x_m)|>2(1-\vert [0.8, 1] \vert )=1.6$.
This contradicts the following fact: 
\begin{eqnarray*}
|\alpha^{m+1}(z)-\alpha(x_m)|&\leq &|\alpha^{m+1}(z)-x_{m+1}|+|x_{m+1}-\widehat\alpha(x_m)|\\
& &\quad +|\widehat\alpha(x_m)-\alpha(x_m)| < 8\varepsilon_1+\delta+\frac{\varepsilon_1}{3}.
\end{eqnarray*}
Thus, $L_\mu\circ c:[0,\nu]\longrightarrow \Sigma$ is a continuous path connecting $L_\mu(L^m(\mathbf{z}))$ with $L_\mu(\mathbf{x}_m)$.
For any $t\in (0,\nu)$,
$$\frac{d}{dt}(L_\mu\circ c)(t)=\left(\frac{\partial \widehat\alpha}{\partial x}(c(t))\dot c_1(t),\ \frac{\partial \beta}{\partial x}(c(t))\dot c_1(t)+\frac{\partial \beta}{\partial y}(c(t))\dot c_2(t)\right).$$
By the condition (\ref{L-condition2}),
$$|[L_\mu(L^m(\mathbf{z}))]_y-[L_\mu(\mathbf{x}_m)]_y|\leq \frac{3}{4\sqrt{2}}\cdot \sqrt{2}\nu <\frac{3\varepsilon}{4},$$
where the `$\sqrt{2}$' of $\sqrt{2}\nu$ is derived from the fact that the maximum of $u+v$ is $\sqrt{2}$ under the assumption of $u^2+v^2=1$.
Note that $[L_\mu(L^m(\mathbf{z}))]_y=[L^{m+1}(\mathbf{z})]_y$, and
$$|[L_\mu (\mathbf{x}_m)]_y-y_{m+1}|\leq |L_\mu (\mathbf{x}_m)-\mathbf{x}_{m+1}|<\delta.$$
It follows that $|[L^{m+1}(\mathbf{z})]_y-y_{m+1}|<7\varepsilon/8$.
Since $|[L^{m+1}(\mathbf{z})]_x-x_{m+1}|=|\alpha^{m+1}(z)-x_{m+1}|<\varepsilon/8$, $|L^{m+1}(\mathbf{z})-\mathbf{x}_{m+1}|<\varepsilon$.
Thus, $\{\mathbf{x}_n\}_{n=1}^\infty$ is $\varepsilon$-shadowed by $\{L^n(\mathbf{z})\}_{n=1}^\infty$.

The proof of (ii) is done similarly by using Proposition \ref{alpha}(ii).
\hspace*{\fill}$\square$

\section{PSSP for Lorenz flows}
In this section, we will prove Theorem \ref{Th2}.
First, consider a Lorenz map $L$ satisfying the conditions (\ref{L-condition2})-(\ref{L-condition4}) and a $L$-Lorenz flow $\varphi$.
Recall that $\Gamma=\{ (0,y,1)\in \Sigma\ ;\ \vert y\vert \leq 1\}$ is the  singularity set on $\Sigma$, and set $\tilde\Gamma=\{ (x,y,z)\in \Pi;\ x=0 \}$, where $\Pi=[-1,1]^2\times [0,1]$.
For the proof of PSSP for $\varphi$, we need to fix a one-parameter family of Lorenz maps $L_\mu$ and $L_\mu$-Lorenz flows.
Here, we suppose that $\{L_\mu\}_{\mu\in [0,\mu_0]}$  is the one-parameter family given in \S 2 and take a $C^1$-one-parameter family $\{\varphi_\mu\}_{\mu\in [0,\mu_0]}$ satisfying $\partial\varphi_\mu/\partial t(\mathbf{x},0)=\partial\varphi/\partial t(\mathbf{x},0)$ for any $\mathbf{x}\in N_{1/10}(\tilde \Gamma,\Pi)$, where $N_\eta (Y,X)$ denotes the $\eta$-neighborhood of a compact subset $Y$ in a metric space $(X,d)$, that is, $N_\eta (Y,X)=\{x\in X\,;\, d(x,Y)\leq \eta\}$.

Let us fix $0<\varepsilon<1$ arbitrarily and determine constants $\hat\delta,\hat\tau >0$ and $\hat\mu\in [0,\mu_0]$ such that any $(\hat\delta,\hat\tau)$-chain for $\varphi_{\hat\mu}$ is $\varepsilon$-shadowed by $\varphi$.

\subsection{Interpolated chains and crossing sequences}
Throughout the remainder of this section, fix $\hat\tau>0$ so that, for any $\mathbf{x} \in \Sigma$ and $\mu\in [0,\mu_0]$, $\varphi_\mu(\mathbf{x},(0,5\hat\tau])\cap \Sigma=\emptyset$.
Let $\{\mathbf{x}_n\}_{n\geq 0}$ be a $(\delta,\hat\tau)$-pseudo-orbit for $\varphi_\mu$, i.e. 
$
\vert
\mathbf{x}_{n+1}-\varphi_\mu(\mathbf{x}_{n},t_n)
\vert\leq \delta
$
for some $\{t_n\}_{n\geq 0}$  with 
$\hat\tau\leq t_n\leq 2\hat\tau$ and $\mathbf{x}_0\in \Sigma$.
Let 
$\{\Phi_{\mu;n}\}_{n\geq 0}$  be the $(\delta,\hat\tau)$-chain for $\varphi_\mu$ associated to $\{\mathbf{x}_n\}_{n\geq 0}$, i.e. 
$\Phi_{\mu;n}=\varphi_\mu(\mathbf{x}_{n}, [0,\ t_n])$.
When $\varphi_\mu(\mathbf{x}_{n}, t_n)\neq \mathbf{x}_{n+1}$, $\sigma_n$ is the open segment in $\mathbb{R}^3$ whose closure connects $\varphi_\mu(\mathbf{x}_{n}, t_n)$ with $\mathbf{x}_{n+1}$, and otherwise $\sigma_n=\emptyset$.
Set
$$
\tilde \Phi_{\mu;n}=\Phi_{\mu;n}\cup \sigma_n.
$$
and call $\{\tilde \Phi_{\mu;n}\}_{n\geq 0}$ the \textit{interpolated $(\delta,\hat\tau)$-chain} for $\varphi_\mu$ associated to $\{\mathbf{x}_n\}_{n\geq 0}$ (or more strictly to $\{\mathbf{x}_n,t_n\}_{n\geq 0}$).
Let $U$ be a small neighborhood of the origin in $\mathbf{R}^3$ with $U\cap \Sigma=\emptyset$.
When $\Phi_{\mu;n}$ is contained in $U$, $\Phi_{\mu;n}$ may have an arbitrarily small length.
On the other hand, $\Phi_{\mu;n}$'s not contained in $U$ have lengths bounded away from zero.
Thus, there exists $\delta_0>0$ such that, for any $\Phi_{\mu;n}$ with $\Phi_{\mu;n}\cap \Sigma\neq \emptyset$, the length of $\Phi_{\mu;n}$ is greater than $3\delta_0$.
This assumption is crucial in our argument below.
In fact, it guarantees that, if $0<\delta\leq \delta_0$, the union $\bigcup_{n\geq 0}\tilde \Phi_{\mu;n}$ of any interpolated $(\delta,\hat\tau)$-chain contains no jagged subsets intersecting $\Sigma$ zigzag.
So, we suppose from now on that $0<\delta\leq \delta_0$.

Let $\{n_i\}_{i\geq 0}$ be the strictly monotone increasing sequence with $n_0=0$ and such that $\{n_i\}_{i\geq 1}$ consists of all positive integers $n$ with $\tilde \Phi_{\mu;n}\cap \Sigma\neq \emptyset$ and $\tilde \Phi_{\mu;n-1}\cap \Sigma= \emptyset$.

\begin{figure}[hbtp]
\centering
\includegraphics[bb=73 530 469 782,clip,width=8.8cm]{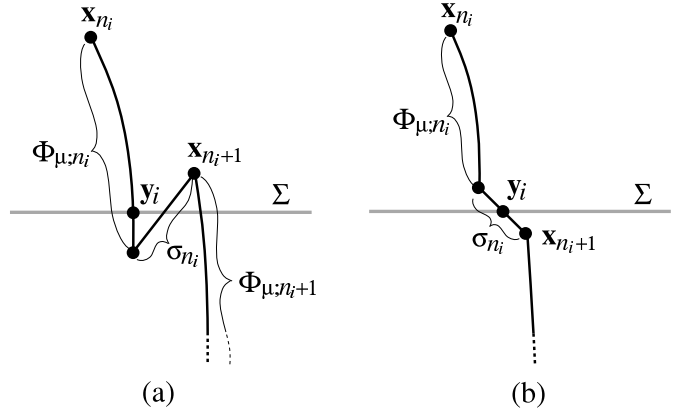}
\caption{In the case (a), both $\Phi_{\mu;n_i}, \Phi_{\mu;n_i+1}$ meet $\Sigma$ non-trivially.
But, the crossing point of $\Phi_{\mu;n_i+1}$ with $\Sigma$ is not an element of $\{\mathbf{y}_i\}_{i\geq 0}$, i.e.\ $n_{i+1}>n_i+1$.}
\label{fig5}
\end{figure}

\begin{definition}\label{crossing}
{\rm
For each $n_i\geq 0$, the \textit{crossing point} $\mathbf{y}_i$ of $\tilde \Phi_{\mu;n_i}$ is a unique point of $\Phi_{\mu;n_i}\cap \Sigma$ if $\Phi_{\mu;n_i}\cap \Sigma\neq \emptyset$, see Fig.\ \ref{fig5}-(a), otherwise $\mathbf{y}_i$ is a point of $\sigma_{n_i}\cap \Sigma$, see Fig.\ \ref{fig5}-(b).
The $\{\mathbf{y}_i\}_{i\geq 0}$ is called the $(\delta,\hat \tau)$-\textit{crossing sequence} for $\varphi_\mu$ associated to $\{\mathbf{x}_n\}_{n\geq 0}$.
}
\end{definition}

\begin{remark}{\rm 
We note that a $(\delta,\hat\tau)$-crossing sequence $\{\mathbf{y}_i\}_{i\geq 0}$ for $\varphi_\mu$ is in general {\it not} a pseudo-orbit for $L_\mu$ even if $\delta>0$ is very small.
The crucial part in our proof of Theorem \ref{Th2} is to show that $\{\mathbf{y}_i\}_{i\geq 0}$ is approximated by a pseudo-orbit $\{\mathbf{w}_i\}_{i\geq 0}$ for $L_\mu$, which in turn is approximated by an actual orbit $\{L^i(\mathbf{z})\}_{i\geq 0}$ of $L$ by Theorem \ref{Th1}.
}
\end{remark}

\subsection{Proof of Theorem \ref{Th2}}
For any $(\delta,\hat \tau)$-crossing sequence $\{\mathbf{y}_{i}\}_{i\geq 0}$ for $\varphi_\mu$, the broken subsegment in $\bigcup_{n\geq 0} \tilde \Phi_{\mu;n}$ connecting 
$\mathbf{y}_i$ with $\mathbf{y}_{i+1}$ is denoted by 
$\left\langle \mathbf{y}_i, \mathbf{y}_{i+1}\right\rangle_\mu^\delta$.
In the case when $\{\mathbf{y}_i\}_{i\geq 0}$ is a finite sequence $\{\mathbf{y}_i\}_{i=0}^m$, $\left\langle \mathbf{y}_m, -\right\rangle^\delta_\mu$ is the broken forward ray in $\bigcup_{n\geq 0} \tilde \Phi_{\mu;n}$ emanating from $\mathbf{y}_m$.

For any $\mathbf{z}\in \Sigma\setminus \Gamma$ and $\mu\in [0,\mu_0]$, let $\tau_{\mu;\mathbf{z}}>0$ be the number with $\varphi_\mu(\mathbf{z},(0,\tau_{\mu;\mathbf{z}}))\cap \Sigma=\emptyset$ and $\varphi_\mu(\mathbf{z},\tau_{\mu;\mathbf{z}})\in \Sigma$, that is, $\varphi_\mu(\mathbf{z},\tau_{\mu;\mathbf{z}})=L_\mu(\mathbf{z})$.

\begin{figure}[hbtp]
\centering
\includegraphics[bb=106 485 487 790,clip,width=8.5cm]{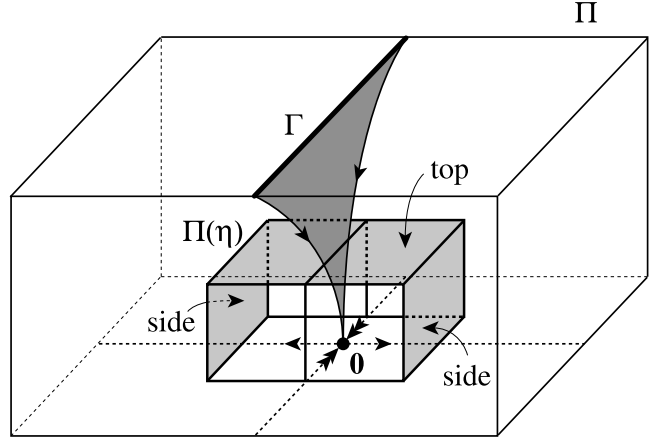}
\caption{`side's represent $\partial_{\rm side}\Pi(\eta)$, and `top' does $\partial_{\rm top}\Pi(\eta)$.
The gray cusp with vertex ${\bf 0}$ is the union $\bigcup_{\mathbf{x}\in \Gamma}\varphi_\mu (\mathbf{x},[0,\infty))\subset W^s({\bf 0})$.}
\label{fig6}
\end{figure}

For any $0<\eta\leq 1$, set $\Pi(\eta)=[-\eta,\eta]^2\times [0,\eta]$, $\partial_{\rm side} \Pi(\eta)=\{-\eta,\eta\}\times [-\eta,\eta]\times [0,\eta]$ and $\partial_{\rm top} \Pi(\eta)=[-\eta,\eta]^2\times \{\eta\}$.
Note that $\partial_{\rm side} \Pi(\eta)$ (resp.\ $\partial_{\rm top} \Pi(\eta)$) consists of two vertical rectangles (resp.\ a single horizontal square) in $\Pi=\Pi(1)$, see Fig.\ \ref{fig6}.
Since any $(\delta,\hat\tau)$-pseudo-orbit $\{\mathbf{x}_n\}_{n\geq 0}$ for $\varphi_\mu$ is taken in the trapping region $T_{\varphi_\mu}$ (see Definition \ref{flowcase}), $\bigcup_{n\geq 0}\tilde \Phi_{\mu;n}\cap \partial \Pi(\eta)$ is contained in $\partial_{\rm top}\Pi(\eta)\cup \partial_{\rm side}\Pi(\eta)$ for any small $\eta>0$, which is suggested by Fig.\ \ref{fig2} and Fig.\ \ref{fig6}.

\begin{figure}[hbtp]
\centering
\includegraphics[bb=144 565 477 782,clip,width=7.4cm]{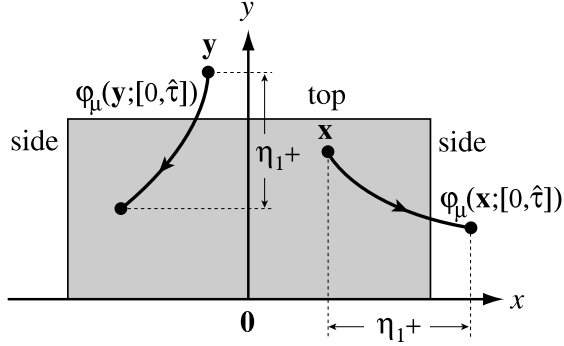}
\caption{The shaded rectangle represents $\Pi(\eta_0)$.}
\label{fig7}
\end{figure}

\begin{figure}[hbtp]
\centering
\includegraphics[bb=119 627 527 787,clip,width=9.1cm]{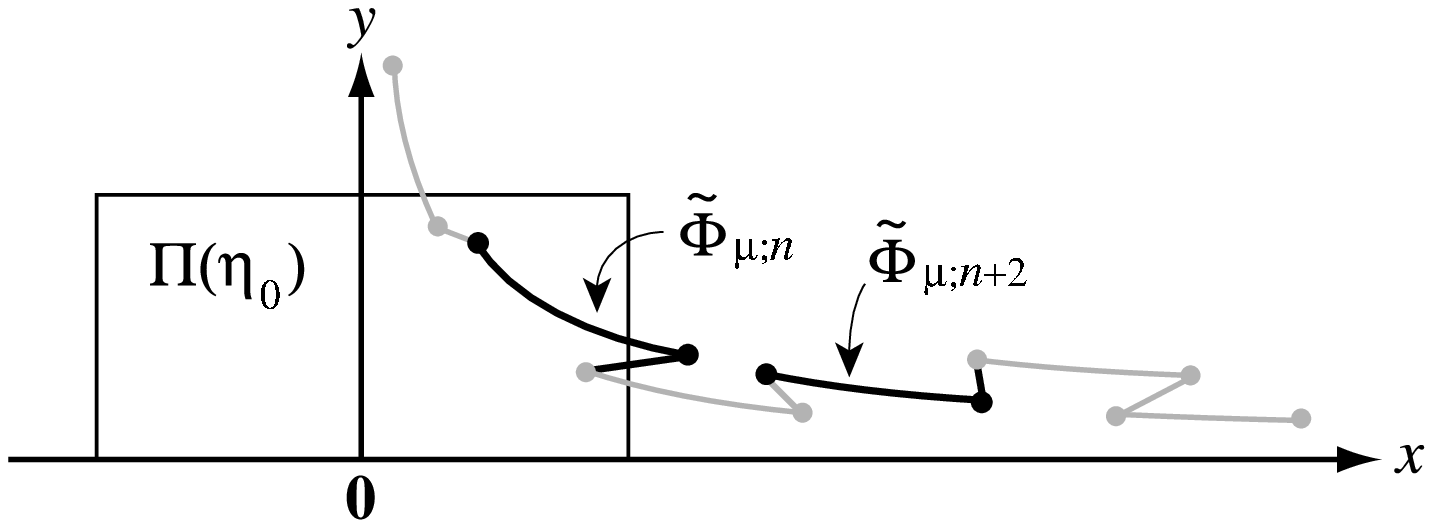}
\caption{}
\label{fig8}
\end{figure}

\begin{figure}[hbtp]
\centering
\includegraphics[bb=67 477 508 780,clip,width=9.8cm]{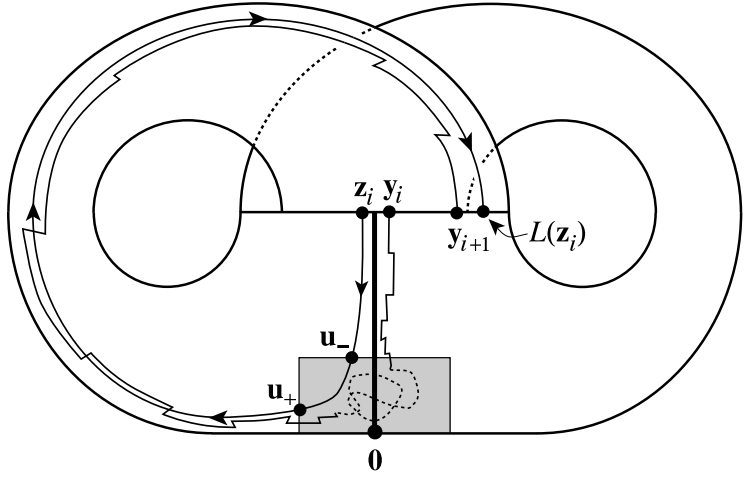}
\caption{The points ${\bf u}_\pm$ represent $\varphi(\mathbf{z}_i,\upsilon_\pm)$.}
\label{fig9}
\end{figure}

\begin{lemma}\label{F1}
There exist $0<\delta_2\leq \delta_0$, $0<\mu_1\leq \mu_0$ and $0<\varepsilon_1<1/10$ such that, for any $0<\mu<\mu_1$ and any $(\delta_2,\hat\tau)$-crossing sequence $\{\mathbf{y}_i\}_{i\geq 0}$ for $\varphi_\mu$, the following conditions {\rm (i)} and {\rm (ii)} hold.
\begin{enumerate}
\renewcommand{\labelenumi}{\rm (\roman{enumi})}
\item
If $\mathbf{z}_i\in \Sigma\setminus \Gamma$ satisfies $\vert \mathbf{z}_i-\mathbf{y}_i\vert \leq {\varepsilon_1}$ and $\vert L(\mathbf{z}_i)-\mathbf{y}_{i+1}\vert \leq {\varepsilon_1}$, then $\langle \mathbf{y}_i,\mathbf{y}_{i+1}\rangle^{\delta_2}_\mu$ is $\varepsilon$-shadowed by $\varphi(\mathbf{z}_i,[0,\tau_{0;\mathbf{z}_i}])$.
\item
If $\{\mathbf{y}_i\}_{i\geq 0}$ is a finite sequence $\{\mathbf{y}_i\}_{i=0}^m$, then for any $\mathbf{z}_m\in \Gamma$ with $\vert \mathbf{z}_m-\mathbf{y}_m\vert \leq \varepsilon_1$, $\langle \mathbf{y}_m,-\rangle^{\delta_2}_\mu$ is $\varepsilon$-shadowed by $\varphi(\mathbf{z}_m,[0,\infty))$.
\end{enumerate}
\end{lemma}
\textit{Proof}.
(i)
From the definition of Lorenz flows, we have $0<\eta_0\leq \varepsilon/30$ such that the restriction $\varphi_\mu\vert{\Pi(3\eta_0)}$ is the linear flow $(e^{\lambda_1 t} x, e^{-\lambda_2 t} y, e^{-\lambda_3 t} z)$ for some $-\lambda_2<-\lambda_3<0<\lambda_3<\lambda_1$ independent of $\mu\in [0,\mu_0]$.
There exists $\eta_1>0$ such that $\vert [\varphi_\mu(\mathbf{x},[0,\hat\tau])]_x\vert-\vert [\mathbf{x}]_x\vert \geq \eta_1$ for any $\mu\in [0,\mu_0]$ and any $\mathbf{x}\in \Pi$ with $\varphi_\mu(\mathbf{x},[0,\hat\tau])\cap \partial_{\rm side}\Pi(\eta_0)\neq \emptyset$ and $[\mathbf{y}]_y -[\varphi_\mu(\mathbf{y},[0,\hat\tau])]_y \geq \eta_1$ for any $\mathbf{y}\in \Pi$ with $\varphi_\mu(\mathbf{y},[0,\hat\tau])\cap \partial_{\rm top}\Pi(\eta_0)\neq \emptyset$, see Fig.\ \ref{fig7}.
Then, there exists $0<\delta_1\leq \min\{\eta_0,\eta_1/2\}$ such that, for any interpolated $(\delta_1,\hat \tau)$-chain $\{\tilde\Phi_{\mu;n}\}_{n\geq 0}$, if $\tilde \Phi_{\mu;n}\cap \partial_{\rm side} \Pi(\eta_0)\neq \emptyset$, then $\tilde \Phi_{\mu;n+2}\cap \Pi(\eta_0)=\emptyset$.
Intuitively, this means that the chain $\{\tilde\Phi_{\mu;k}\}_{k\geq n}$ eventually goes away from $\Pi(\eta_0)$ if $\tilde \Phi_{\mu;n}\cap \partial_{\rm side} \Pi(\eta_0)\neq \emptyset$, see Fig.\ \ref{fig8}.
Similarly, one can choose the $\delta_1$ so that, if $\tilde \Phi_{\mu;n}\cap \partial_{\rm top}\Pi(\eta_0)\neq \emptyset$, then $\tilde\Phi_{\mu;n+2}\cap \partial_{\rm top}\Pi(\eta_{0})=\emptyset$.
For any $\mathbf{z}\in \Sigma\setminus \Gamma$, we set $t_{\mu;\mathbf{z}}=0$ if $\varphi_\mu(\mathbf{z},[0,\tau_{\mu;\mathbf{z}}])\cap \Pi(\eta_0)=\emptyset$ and otherwise $t_{\mu;\mathbf{z}}=\tau_+-\tau_-$, where $[\tau_-,\tau_+]$ is the subinterval of $[0,\tau_{\mu;\mathbf{z}}]$ with $\varphi_\mu(\mathbf{z},[\tau_-,\tau_+])=\varphi_\mu(\mathbf{z},[0,\tau_{\mu;\mathbf{z}}])\cap \Pi(\eta_0)$.
Since $\varphi_\mu$ has no singular points in $T_\mu\setminus \Pi(\eta_0)$, there exists $s_0>0$ such that, for any $\mathbf{z}\in \Sigma\setminus \Gamma$ and $\mu\in [0,\mu_0]$, $\tau_{\mu;\mathbf{z}}-t_{\mu;\mathbf{z}}<s_0$.
From this, we know that, for any interpolated $(\delta_1,\hat\tau)$-chain $\{\Phi_{\mu;n}\}_{n\geq 0}$, 
there exists the number of $n$'s with $n_i\leq n\leq n_{i+1}$ 
such that $\tilde\Phi_{\mu;n}$ 
is not wholly contained in $\Pi(\eta_0)$,  
and  is bounded by a constant independent of $\mathbf{z}\in \Sigma\setminus \Gamma$ and $\mu\in [0,\mu_0]$.
Then, one can choose $0<\delta_2\leq \delta_1$, $0<\mu_1\leq \mu_0$, $0<\varepsilon_1<1$ such that, for any $(\delta_2,\hat \tau)$-crossing sequence $\{\mathbf{y}_i\}_{i\geq 0}$ for $\varphi_\mu$, $\langle \mathbf{y}_i,\mathbf{y}_{i+1}\rangle_\mu^{\delta_2}\setminus \Pi(\eta_0)$ is $\varepsilon$-shadowed by $\varphi(\mathbf{z}_i,[0,\tau_{0;\mathbf{z}_i}]\setminus [\upsilon_-,\upsilon_+])$ if $0<\mu\leq \mu_1$, $\vert \mathbf{z}_i-\mathbf{y}_i\vert \leq {\varepsilon_1}$ and $\vert L(\mathbf{z}_i)-\mathbf{y}_{i+1}\vert \leq {\varepsilon_1}$, where $[\upsilon_-,\upsilon_+]$ is the subinterval (possibly empty) of $[0,\tau_{0;\mathbf{z}_i}]$ with $\varphi(\mathbf{z}_i,[0,\tau_{0;\mathbf{z}_i}])\cap \Pi(\eta_0)=\varphi(\mathbf{z}_i,[\upsilon_-,\upsilon_+])$.
Since the diameter of $\Pi(\eta_0)$ is less than $\varepsilon/2$, $\langle \mathbf{y}_i,\mathbf{y}_{i+1}\rangle_\mu^{\delta_2}\cap \Pi(\eta_0)$ is also $\varepsilon$-shadowed by $\varphi(\mathbf{z}_i,[\upsilon_-,\upsilon_+])$, see Fig.\ \ref{fig9}.
This shows the assertion (i).\\
(ii)
The proof is quite similar to that of (i).
Suppose that $0<\mu\leq \mu_1$ and $\{\mathbf{y}_i\}_{i=0}^m$ is a finite $(\delta_2,\hat\tau)$-crossing sequence for $\varphi_\mu$.
From the argument in (i), $\langle \mathbf{y}_m,-\rangle_\mu^{\delta_2}$ is disjoint from $\partial_{\rm side}\Pi(\eta_0)$.
For any $\mathbf{z}_m\in \Gamma$ with $\vert \mathbf{z}_m-\mathbf{y}_m\vert \leq \varepsilon_1$, let $\upsilon_-$ be a unique point in $[0,\infty)$ with $\varphi(\mathbf{z}_m,\upsilon_-)\in \partial_{\rm top}\Pi(\eta_0)$.
Then, $\langle \mathbf{y}_m,-\rangle^\mu_{\delta_2}\setminus\Pi(\eta_0)$ is $\varepsilon$-shadowed by $\varphi(\mathbf{z},[0,\upsilon_-])$.
Since $\varphi(\mathbf{z}_i,[\upsilon_-,\infty))$ is contained in $\Pi(\eta_0)\cap \tilde \Gamma$, $\langle \mathbf{y}_m,-\rangle_\mu^{\delta_2}\cap \Pi(\eta_0)$ is $\varepsilon$-shadowed by $\varphi(\mathbf{z}_m,[\upsilon_-,\infty))$.
Thus, $\langle \mathbf{y}_m,-\rangle_\mu^{\delta_2}$ is $\varepsilon$-shadowed by $\varphi(\mathbf{z}_m,[0,\infty))$.
\hspace*{\fill}$\square$\\

By Theorem \ref{Th1}, there exist $\hat \mu\in (0,\mu_1]$ and $\xi_0>0$ such that any $\xi_0$-pseudo-orbit for $L_{\hat\mu}$ is $\varepsilon_1/2$-shadowed by an actual orbit for $L$.
From now on, we fix  a $\hat\mu>0$ satisfying this condition and suppose that any pseudo-orbits and crossing sequences are those for $\varphi_{\hat\mu}$.
Here, one can suppose that the $\xi_0$ is less than $\varepsilon_1$.\\

\noindent\textit{Proof of Theorem \ref{Th2}.}
First, let us consider the case when crossing sequences associated with a pseudo-orbits are infinite.
We will show that there exists $0<\hat\delta\leq \delta_2$ such that, for the infinite crossing sequence $\{\mathbf{y}_i\}_{i\geq 0}$ associated with a $(\hat\delta,\hat\tau)$-pseudo-orbit $\{\mathbf{x}_n\}_{n\geq 0}$, there is an infinite sequence $\{\mathbf{w}_i\}_{i\geq 0}$ in $\Sigma$ which is a $\xi_0$-pseudo-orbit for $L_{\hat\mu}$ satisfying
\begin{equation}\label{quoter}
\vert \mathbf{y}_i-\mathbf{w}_i\vert <\varepsilon_1/2
\end{equation}
for any $i\geq 0$.

Note that any flow of $\varphi_{\hat\mu}$ emanating from $\mathbf{0}$ tends toward either $\mathbf{v}_+$ or $\mathbf{v}_-$.
Take $0<\eta_2\leq \eta_0$ such that, for any $\mathbf{z}\in \Pi(\eta_2)\setminus \tilde \Gamma$, the first crossing point of $\varphi_{\hat\mu}(\mathbf{z},t); t> 0$ with $\Sigma$ is contained in either $N_{\xi_0/3}(\mathbf{v}_+,\Sigma)$ or $N_{\xi_0/3}(\mathbf{v}_-,\Sigma)$.
There exists $0<\delta_3\leq \delta_2$ such that, for any $(\delta_3,\hat\tau)$-crossing sequence $\{\mathbf{y}_i\}_{i\geq 0}$, if $\langle \mathbf{y}_i,\mathbf{y}_{i+1}\rangle_{\hat\mu}^{\delta_3}\cap \Pi(\eta_2)\neq \emptyset$, then $\mathbf{y}_{i+1}$ is contained in either $N_{\xi_0/2}(\mathbf{v}_+,\Sigma)$ or $N_{\xi_0/2}(\mathbf{v}-,\Sigma)$.
Then, we have $0<\xi_1\leq \xi_0/4$ and $0<\delta_4\leq \delta_3$ such that, for any $(\delta_4,\hat\tau)$-crossing sequence $\{\mathbf{y}_i\}_{i\geq 0}$, $\langle \mathbf{y}_i,\mathbf{y}_{i+1}\rangle_{\hat\mu}^{\delta_4}$ meets $\Pi(\eta_2)$ non-trivially if $\vert  [\mathbf{y}_i]_x\vert \leq \xi_1$.
One can take $0<\hat \delta\leq \delta_4$ such that if $\vert [\mathbf{y}_i]_x\vert\geq \xi_1$ for a $(\hat\delta,\hat\tau)$-crossing sequence $\{\mathbf{y}_i\}_{i\geq 0}$, then $\langle \mathbf{y}_i,\mathbf{y}_{i+1}\rangle_{\hat\mu}^{\hat\delta}$ is disjoint from $\tilde \Gamma$ and
\begin{equation}\label{pseudo}
\vert \mathbf{y}_{i+1}-L_{\hat\mu}(\mathbf{y}_i)\vert <\xi_0/2.
\end{equation}
We set $\mathbf{y}_i=\mathbf{w}_i$ if $\vert [\mathbf{y}_i]_x\vert\geq \xi_1$.
Here, we need to consider the following three cases.

\begin{description}
  \item[Case 1] $\vert [\mathbf{y}_i]_x\vert\geq \xi_1$ and $\vert [\mathbf{y}_{i+1}]_x\vert\geq \xi_1$.

Since $\mathbf{y}_i=\mathbf{w}_i$ and $\mathbf{y}_{i+1}=\mathbf{w}_{i+1}$, by (\ref{pseudo}), $\vert \mathbf{w}_{i+1}-L_{\hat\mu}(\mathbf{w}_i)\vert <\xi_0$.
  \item[Case 2] $\vert [\mathbf{y}_i]_x\vert\leq \xi_1$.

In this case, $\langle \mathbf{y}_i,\mathbf{y}_{i+1}\rangle_{\hat\mu}^{\hat\delta}$ may intersect with $\tilde \Gamma$ non-trivially.
Then, it can happen that $\mathbf{y}_{i+1}\in N_{\xi_0/2}(\mathbf{v}_{\iota},\Sigma)$ and $L_{\hat\mu}(\mathbf{y}_i)\in N_{\xi_0/2}(\mathbf{v}_{-\iota},\Sigma)$ for some $\iota\in \{+,-\}$, see Fig.\ \ref{fig10}.
\begin{figure}[h]
\centering
\includegraphics[bb=91 368 485 815,clip,width=8.8cm]{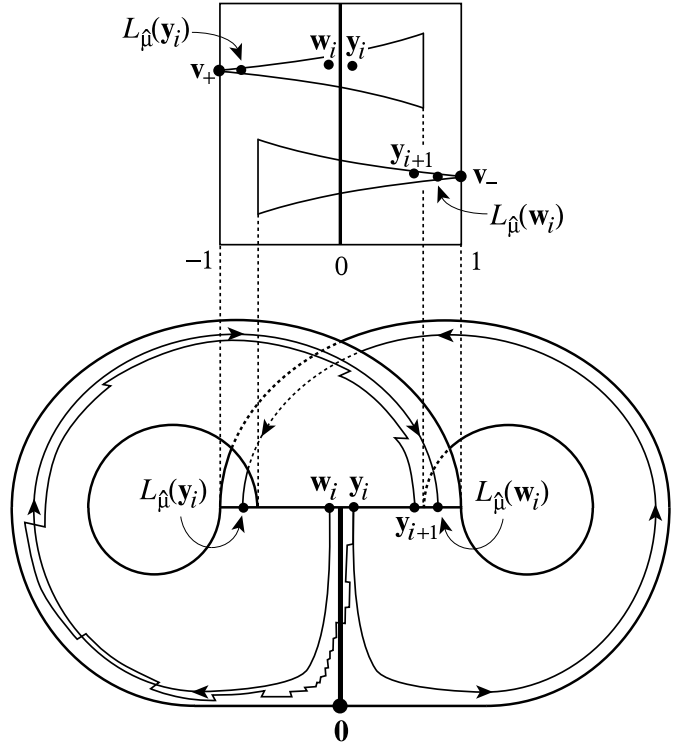}
\caption{The case of $[\mathbf{y}_i]_x>0$, $[\mathbf{y}_{i+1}]_x>0$ and $[\mathbf{w}_i]_x<0$.
Then, $L_{\hat\mu}(\mathbf{y}_i)$ is not approximated by $\mathbf{y}_{i+1}$.}
\label{fig10}
\end{figure}
Take a point $\mathbf{w}_i\in \Sigma$ with $[\mathbf{w}_i]_y=[\mathbf{y}_i]_y$, $0<\vert [\mathbf{w}_i]_x\vert\leq \xi_1$ and $\iota={\rm sign}[\mathbf{w}_i]_x=-{\rm sign}[\mathbf{y}_{i+1}]_x$.
This definition implies $\vert \mathbf{w}_i-\mathbf{y}_i\vert \leq 2\xi_1\leq \xi_0/2<\varepsilon_1/2$.
Since $\vert [\mathbf{y}_{i+1}]_x\vert \geq 1-\xi_0/2>\xi_1$, $\mathbf{y}_{i+1}=\mathbf{w}_{i+1}$ and hence
\begin{eqnarray*}
\vert L_{\hat\mu}(\mathbf{w}_i)-\mathbf{w}_{i+1}\vert&=&\vert L_{\hat\mu}(\mathbf{w}_i)-\mathbf{y}_{i+1}\vert\\
&\leq& \vert L_{\hat\mu}(\mathbf{w}_i)-\mathbf{v}_{\iota}\vert +\vert \mathbf{v}_{\iota}-\mathbf{y}_{i+1}\vert \\
&<& \xi_0/2+\xi_0/2=\xi_0.
\end{eqnarray*}
  \item[ Case 3] $\vert [\mathbf{y}_i]_x\vert\geq \xi_1$ and $\vert [\mathbf{y}_{i+1}]_x\vert< \xi_1$.

As was shown in the argument of Case 2, $\vert \mathbf{w}_{i+1}-\mathbf{y}_{i+1}\vert \leq 2\xi_1$.
Since $\mathbf{y}_i=\mathbf{w}_i$, the inequality (\ref{pseudo}) implies
\begin{eqnarray*}
\vert L_{\hat\mu}(\mathbf{w}_i)-\mathbf{w}_{i+1}\vert&\leq &\vert L_{\hat\mu}(\mathbf{y}_i)-\mathbf{y}_{i+1}\vert +\vert \mathbf{y}_{i+1}-\mathbf{w}_{i+1}\vert\\
&<& \xi_0/2+2\xi_1\leq \xi_0.
\end{eqnarray*}
\end{description}

By Cases 1-3, $\{\mathbf{w}_i\}_{i\geq 0}$ is a $\xi_0$-pseudo-orbit of $L_{\hat\mu}$ satisfying (\ref{quoter}).
By Theorem \ref{Th1}\,(i), there exists $\mathbf{z}\in \Sigma\setminus \Gamma$ with $\bigcup_{i=0}^\infty L^i (\mathbf{z})\cap \Gamma =\emptyset$ and such that $\{L^i(\mathbf{z})\}_{i\geq 0}$ $\varepsilon_1/2$-shadows $\{\mathbf{w}_i\}_{i\geq 0}$.
Since $\vert \mathbf{y}_i-L^i(\mathbf{z})\vert \leq \vert \mathbf{y}_i-\mathbf{w}_i\vert +\vert \mathbf{w}_i-L^i(\mathbf{z})\vert <\varepsilon_1$, by Lemma \ref{F1}\,(i), the $(\hat\delta,\hat\tau)$-chain $\{\Phi_{\hat\mu;n}\}_{n\geq 0}$ associated to $\{\mathbf{x}_n\}_{n\geq 0}$ is $\varepsilon$-shadowed by the actual orbit $\varphi(\mathbf{z},t);t\geq 0$.

Next, we consider the case when crossing sequences associated with $(\delta_2,\hat\tau)$-pseudo-orbits is finite.
By the argument as above, we have $0<\hat\delta\leq \delta_2$ such that, for the finite crossing sequence $\{\mathbf{y}_i\}_{i=0}^m$ associated with any $(\hat\delta,\hat\tau)$-pseudo-orbit $\{\mathbf{x}_n\}_{n\geq 0}$, there is a sequence $\{\mathbf{w}_i\}_{i=0}^m$ in $\Sigma$ which is a $\xi_0$-pseudo-orbit for $L_{\hat\mu}$ satisfying $\vert \mathbf{y}_i-\mathbf{w}_i\vert \leq \varepsilon_1/2$ for any $i\in \{0,1,\cdots,m\}$.
By Theorem \ref{Th1}\,(ii), there exists $\mathbf{z}\in \Sigma$ with $\bigcup_{i=0}^{m-1}L^i(\mathbf{z})\cap \Gamma=\emptyset$, $L^m(\mathbf{z})\in \Gamma$ and such that $\{L^i(\mathbf{z})\}_{i=0}^m$ $\varepsilon_1/2$-shadows $\{\mathbf{w}_i\}_{i=0}^m$.
Then, by applying Lemma \ref{F1}\,(i) $(m-1)$-times and (ii) once, one can show that the $(\hat\delta,\hat\tau)$-chain $\{\Phi_{\hat\mu;n}\}_{n\geq 0}$ associated to $\{\mathbf{x}_n\}_{n\geq 0}$ is $\varepsilon$-shadowed by the actual orbit $\varphi(\mathbf{z},t);t\geq 0$.
This completes the proof of Theorem \ref{Th2}.
\hspace*{\fill}$\square$\\

\bibliographystyle{amsplain}

\end{document}